\documentclass[12pt]{article}
\usepackage{amssymb,amsmath,amscd,amsthm}
\usepackage{graphicx}

\newtheorem{theorem}{Theorem}[section]

\newtheorem{lemma}[theorem]{Lemma}

\date{}

\begin{document}

\title{Spectral analysis of non-local Schr\"{o}dinger operators. }

\author{Yu. Kondratiev \footnote{Fakultat fur Mathematik, Universitat Bielefeld, 33615 Bielefeld, Germany,
kondrat@math.uni-bielefeld.de. The work was partially supported  by the DFG through
SFB 701, Bielefeld University.}, S. Molchanov \footnote{Dept of Mathematics and Statistics, UNC at Charlotte, NC 28223 and National Research Univ., Higher School of Economics, Russian Federation, smolchan@uncc.edu. The work was partially supported  by the NSF grant DMS-1410547 and by the DFG through
SFB 701, Bielefeld University.
}, B. Vainberg \footnote{Dept of Mathematics and Statistics, UNC at Charlotte, NC 28223, brvainbe@uncc.edu. The work was partially supported  by the NSF grant DMS-1410547; corresponding author.}}

\maketitle

\begin{abstract} We study spectral properties of convolution operators $\mathcal L$ and their perturbations $H=\mathcal L+v(x)$ by compactly supported potentials. Results are applied to determine the front propagation of a population density governed by operator $H$ with a compactly supported initial density provided that $H$ has positive eigenvalues. If there is no positive spectrum, then the stabilization of the population density is proved.
\end{abstract}



 \section {Introduction}
In this paper we study the spectral properties of the Hamiltonian $H$ in $L^2(R^d)$ given by
\begin{equation}\label{HHH}
H=\mathcal L+v(x), \quad \mathcal L\psi(x)=\chi\int_{y\in R^d}(\psi(x+y)-\psi(x))a(y)dy, \quad x\in R^d,
\end{equation}
where
\begin{equation}\label{aone1}
a(y)=a(-y), \quad a\geq 0, \quad \int_{R^d}a(y)dy=1,
\end{equation}
and $v(x)\geq 0$ is continuous and compactly supported.

The operator $\mathcal L$ is the generator of a symmetric random walk on $R^d$ with the intensity of jumps equal to $\chi>0$. Function $a(y)$ is the density of transition from $x$ to $x+y$ at the moment of the jump. Operators $H$ of the form (\ref{HHH}) appear in many applications, such as models of population dynamics that include the KPP type processes (where the offspring start at the location of the parent particle \cite{kpp}, \cite{d1}, \cite{d2}, \cite{ger}) and contact processes (where the locations of the new particles are
randomly distributed around the location of the parent \cite{kkp}, \cite{kosk}). Another area of application is the theory of phase transitions for homopolymers, \cite{CrkoMoV}, \cite{CrkoMoV2}. The central object of investigation in both applications is the solution $u(t,x)$ of the parabolic problem
\begin{equation}\label{xx}
\frac{\partial u}{\partial t}=Hu=(\mathcal L+v(x))u, \quad u(0,x)=u_0(x).
\end{equation}

In the population dynamics models, $u(t,x)$ is the first correlation function, i.e., the density of the population. The potential $v(x)$ in these models is the difference between the birth and death rates: $v(x)=\beta(x)-\mu(x)$. In the KPP type models, $\chi, \beta$ and $\mu$ are unrelated, while in the contact models, $\beta=\chi=$const., $\mu\leq \beta$ (see \cite{kosk}). In the theory of homopolymers, the solution of (\ref{xx}) is the partition function of the Gibbs distribution.

If $v(x)\equiv 0$ (i.e., the birth and death rates are equal), then the homogeneous equation (\ref{xx}) has a solution $u(t,x)\equiv $const.
In the presence of a non-negative potential, the problem
\begin{equation}\label{u0c}
\frac{\partial u}{\partial t}=(\mathcal L+v(x))u, \quad \quad u(0,x)=C={\rm const.},
\end{equation}
may exhibit two significantly different types of large time behavior. If operator $H$ has a positive eigenvalue $\lambda_0$ (the largest eigenvalue has a positive eigenfunction, ground state, $\phi_0(x))$, then the solution of the problem (\ref{u0c}) grows exponentially in time on any fixed bounded domain $D\subset R^d$. This is a manifestation of an instability under local perturbation (by the potential $v$). If the spectrum of $H$ belongs to $(-\infty,0]$, then in many cases the density $u(t,x)$ remains bounded as $t\to\infty$, and one can expect the existence of a steady state for the perturbed problem.

The fundamental difference between the classical Schr\"{o}dinger operator and the Hamiltonian $H$ under consideration is that the potential term in the Schr\"{o}dinger operator is a relatively compact perturbation of the Laplacian, but now this term is not relatively compact with respect to the operator $\mathcal L$. In order to overcome this difficulty, an essential restriction on $v$ will be imposed: in majority of cases we will assume that $0\leq v(x)<\chi$. In the context of population models, this condition
precludes the scenario where the growth of the population in a single point overwhelms the ability of the
population to spread. It automatically holds in the case of the contact models.

The problem of the existence of a ground state energy $\lambda_0(v)>0$ was discussed in the recent paper \cite{kpz}. The central idea of the approach in  \cite{kpz} is similar to ideas used in \cite{CrkoMoV}, \cite{CrkoMoV2}, but the setting is different. Let us stress that the results in the paper \cite{kpz} were based solely on the Perron-Frobenius theory. Here we impose additionally the symmetry condition $a(y)=a(-y)$ that leads to the self-adjointness of $\mathcal L$ and allows us to provide a detailed spectral analysis of $\mathcal L$ or $H$.

In this paper we will study:

a) The a.c. spectrum of $\mathcal L$ and $H$. Conditions for the spectrum to be pure a.c..

b) Examples of operators  $\mathcal L$  with the point spectrum imbedded into the a.c. one.

c) The discrete spectrum outside of the a.c. one; in particular, the difference in the properties of this spectrum for operators $H$ with recurrent and transient underlying random walks (an analogue of the properties that hold for the standard Schr\"{o}dinger operators).

d) The asymptotic behavior at infinity of the ground state $\phi_0(x)$.

e) The propagation of the population front when the ground state exists and convergence of the solution of (\ref{u0c}) to a bounded solution as $t\to\infty$ in the case of absence of positive spectrum.

Certain results for operator $H$ on the lattice $Z^d$ (where the situation is similar and simpler) can be found
 in \cite{MW}, \cite{yar}, \cite{Myar}.

 \section {Spectral analysis of operator $\mathcal L$}

After the Fourier transform, operator $\mathcal L$ becomes the operator of multiplication by function $\chi(\widehat{a}(k)-1)$, where $\widehat{a}=\widehat{a}(k)$ is given by
 \begin{equation}\label{ahat}
\widehat{a}(k)=\int_{R^d}e^{-i(k,y)}a(y)dy=\int_{R^d}\cos(k,y)a(y)dy.
\end{equation}
From (\ref{aone1}) it follows that function $\widehat{a}$ is real-valued and
\begin{equation}\label{cond1}
\widehat{a}(0)=1;~~~|\widehat{a}(k)|< 1, ~~k\neq 0; \quad  \widehat{a}(k)\to 0 ~~{\rm as}~~k\to\infty.
\end{equation}
We will assume that function $\widehat{a}$ satisfies the following condition:
\begin{equation}\label{cond}
\widehat{a}\in L^1(R^d),
\end{equation}
which is natural to guarantee  the continuity of the jump distribution density $a(y)$.

Consider a more general operator than $\widehat{\mathcal L}$. Let $B:L^2(R^d)\to L^2(R^d)$ be an operator of multiplication by a smooth real-valued function $\beta=\beta(k) $, i.e., $B\phi(k)=\beta (k)\phi(k).$
Denote the $\lambda$-level set of the function $\beta$ by $\beta_\lambda$, i.e.,
\begin{equation}\label{all}
\beta_\lambda=\{k\in R^d: \beta(k)=\lambda\}.
\end{equation}If $S$ is a set in $R^d$, then its Lebesgue measure will be denoted by $m(S)$, and $m_1(S)$ will be used instead of $m(S)$ when $d=1$ (notation $\mu$ is preserved for the spectral measure of operator $B$).
\begin{lemma}\label{laaa}
Operator $B$ has the following properties.

(1) The spectrum of $B$ coincides with the closure of the range of function $\beta$.

(2) The point spectrum of $B$ consists of points $\lambda$ for which $m(\beta_\lambda)>0$. All the eigenvalues have infinite multiplicity.

(3) If for some open interval $\Delta\subset {\rm Sp}B$,
\begin{equation}\label{333}
m(S_\Delta)=0 \quad  where \quad S_\Delta=\{k\in R^d:  \beta(k)\in\Delta,\nabla\beta(k)=0\},
\end{equation}
then ${\rm Sp}B$ is absolutely continuous on $\Delta$.

4) If $\beta$ is analytic and not constat, then the spectrum of $B$ is absolutely continuous.
\end{lemma}
{\bf Proof.}\footnote {The proof uses P. Kuchment's ideas suggested in our conversation.} The first statement is obvious. The second one is also trivial. Indeed, each function with the support in the set $\beta_\lambda$ is an eigenfunction of $B$ with the eigenvalue $\lambda$. Conversely, if $(\beta(k)-\lambda)f(k)=0,$ then $f(k)=0$ for $k\notin \beta_\lambda$, and therefore $f=0$ as an element of $L^2(R^d)$ if $m(\beta_\lambda)=0$. Let us prove the third statement. One needs only  to show that the spectral measure $\mu_f$ of an arbitrary element $f\in L^2(R^d)$ does not have a singular continuous component in $\Delta$.

Let $P_\delta$ be the spectral projection of operator $B$ on an arbitrary interval $\delta\subset\Delta$. Then $\mu_f(\delta)=<P_\delta f,f>$.
Since $\Delta$ does not contain eigenvalues, by the Stone formula, the following relation holds for each smooth compactly supported function $f$:
\[
P_\delta f(k)=\lim_{\varepsilon\to+0}\frac{1}{\pi}\int_\delta\frac{\varepsilon f(k)d\lambda}{(\beta(k)-\lambda)^2+\varepsilon^2}=\chi_\delta(k)f(k),
\]
where $\chi_\delta$ is the indicator function of the set $\{k:\beta(k)\in \delta\}$. Thus the same relation is valid for arbitrary $f\in L^2(R^d)$. Hence
\[
\mu_f(\delta)=\int_{R^d}\chi_\delta(k)|f(k)|^2dk,
\]
and therefore, for each Borel-measurable set $\gamma\subset\Delta,$
\[
\mu_f(\gamma)=\int_{R^d}\chi_\gamma(k)|f(k)|^2dk, \quad \chi_\gamma(k)=      1 ~{\rm if}~ \beta(k)\in \gamma, ~~\chi_\gamma(k)=
                                                                             0  ~{\rm if} ~ \beta(k)\notin \gamma.
\]
It remains to show that $\mu_f(\gamma)=0$ if $m_1(\gamma)=0$.

Assume that there exists a set $\gamma\subset\Delta$ of Lebesgue measure zero such that
\[
\mu_f(\gamma)=\int_{R^d}\chi_\gamma(k)|f(k)|^2dk=\varepsilon>0.
\]
Using an approximation of $f$ in $L^2(R^d)$ by bounded functions with compact supports, one can find $r<\infty$ and $\varepsilon_1>0$ such that
\begin{equation}\label{mmm}
\int_{|k|<r}\chi_\gamma(k)dk\geq\varepsilon_1>0.
\end{equation}
We split the support of $\chi_\gamma$ in the ball $|k|<r$ into two parts: $b_1(\sigma)$ is the subset of the support of $\chi_\gamma$ where $|\nabla \beta(k)|\geq\sigma>0$ and $b_2(\sigma)$ is the subset where $|\nabla \beta(k)|<\sigma$. Obviously $m(b_1(\sigma))=0$ and $m(b_2(\sigma))\to 0$ as $\sigma\to 0$. Thus the support of $\chi_\gamma$ in the ball $|k|<r$ has Lebesgue measure zero. This contradicts (\ref{mmm}) and completes the proof of the third statement of the lemma.

The last statement follows from the previous one and the Weierstrass preparation theorem. Indeed, let $\nabla \beta(k_0)=0$. We choose $i$ such that $\frac{\partial\beta}{\partial k_i}$ does not vanish identically in a neighborhood of $k=k_0$.  Then there is a unitary transformation $U$ in $R^d$ such that the function $\frac{\partial\beta}{\partial k_i}$  does not vanish identically on the $z_1$-axis, where $z=U(k-k_0)$. Hence the Weierstrass preparation theorem implies that the equation $\frac{\partial\beta}{\partial k_i}=0$  is equivalent to  $W(z)=0$ in a small neighborhood of $z=0$, where $W(z)=z_1^N+g_1(z')z_1^{N-1}+...+g_N(z')$ is a Weierstrass polynomial in $z_1$ whose coefficients are analytic functions with respect to the remaining variables $z'\in R^{d-1}$. For each fixed value of $z'$, the polynomial $W(z)$ has $N$ roots in $z_1$. Thus the set $\{k:~\frac{\partial\beta}{\partial k_i}=0,|k-k_0|<\delta\}$, where $\delta>0$ is small enough, has zero measure. This implies (\ref{333}) and allows one to apply statement 3 of the lemma.

\qed

The following three classes of densities $a(y)$ will be considered.

1. {\it Processes with light tails} for which $|a(y)|<Ce^{-\delta|y|}$, and therefore $\widehat{a}=\widehat{a}(k)$ is analytic in $k$ when $|{\rm Im}k|<\delta$ with some $\delta>0$.

2. {\it Moderate tails} are defined by the condition

\begin{equation}\label{case1}
a(y)=\frac{a_0(\dot{y})}{|y|^{d+\gamma}}(1+O(|y|^{-\epsilon})), \quad  y\to\infty, \quad \dot{y}=\frac{y}{|y|},\quad \gamma>2,
\end{equation}
where $a_0$ is assumed to be continuous and positive.

In both cases of the light tails and the moderate tails, the density $a(y)$
has second moments, and therefore $\widehat{a}\in C^2$. We assume that the the covariance matrix $B$ is not degenerate:
\begin{equation}\label{detb}
\det B\neq 0,  \quad B=[-\frac{\partial^2\widehat{a}(k)}{\partial k_i\partial k_j}]|_{k=0}.
\end{equation}

3. {\it Processes with heavy tails} are defined by (\ref{case1}) with $\gamma\in(0,2)$. Thus the second moments do not exist in this case. We will assume that $\widehat{a}(k)\in C^1$ when $k\neq 0$. This assumption holds \cite{amv} if $\gamma>1$ or $0<\gamma\leq 1,~a_0(\dot{y})$ is smooth and a two term asymptotic expansion is valid instead of (\ref{case1}).

Relations (\ref{cond1}) imply that the closure of the range of the function $\chi(\widehat{a}(k)-1)$ is $[-a,0]$, where $\chi\leq a\leq 2\chi$.
\begin{theorem}\label{ttt}
The spectrum of operator $\mathcal L$ coincides with the segment $[-a,0]$. The spectrum of $\mathcal L$ is pure a.c. in the case of light tails. It is a.c. in a neighborhood of the origin in the case of moderate tails. It may contain a countable set of embedded eigenvalues of infinite multiplicity in the case of moderate and heavy tails.
\end{theorem}
{\bf Remark.} An example with embedded eigenvalues will be constructed when $\widehat{a}\in C_0^\infty$, and therefore, $a(y)$ is analytic and decays at infinity faster than any power.

{\bf Proof.} The first two statements are proved in Lemma \ref{laaa}. Let us prove that the spectrum is a.c. in a neighborhood of the origin in the case of moderate tails. Note that $\widehat{a}\in C^2$. Since $\widehat{a}(0)=1$ and $\nabla\widehat{a}(0)=0$ (due to the symmetry of $a(y)$), we have $\widehat{a}(k)=1-\frac{1}{2}(Bk,k)+o(|k|^2),~k\to 0,$ where the covariance matrix $B$ is non-degenerate. Therefore $\nabla\widehat{a}(k)=-Bk+o(|k|),~k\to 0,$ and function $\widehat{a}$ does not have critical points in a $\delta$-neighborhood of the origin other than $k=0$.

The strict inequality $\widehat{a}(k)<1$ holds if $k\neq 0$, since otherwise
\[
\int_{R^d}(1-\cos(k_0,y))a(y)dy=0
\]
for some $k_0\neq 0$. The latter implies $\cos(k_0,y)=1$ on an open set where $a(y)>0$. This contradicts the analyticity of $\cos(k_0,y)$. Since $\widehat{a}(k)<1$ for $k\neq 0$ and $\widehat{a}(k)\to 0$ as $k\to \infty$, there exists $\varepsilon >0$ such that $\widehat{a}(k)<1-\varepsilon$ when $|k|>\delta$. Thus item (3) of Lemma \ref{laaa} can be applied to the interval $\Delta=(-\chi\varepsilon,0)$, i.e., the spectral measure is a.c. there.

In order to complete the proof of the theorem, we need to construct an example of an operator with embedded eigenvalues.
\begin{figure}[htbp]
\begin{center}
\includegraphics[width=0.4\columnwidth]{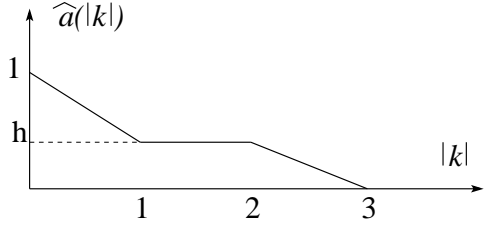}
\end{center}
\caption{Function $\widehat a(|k|)$.}
\label{fig-1}
\end{figure}
 First we will construct a one-dimensional operator with at least one eigenvalue of infinite multiplicity. Consider the function $\widehat{a}=\widehat{a}(|k|),~k\in R^1,$ whose graph is given in Fig 1. We will make it smoother later. The inverse Fourier transform gives
\[
a(x)=\frac{1}{\pi x^2}[(1-\cos x)(1-h)+(\cos 2x-\cos 3x)h].
\]
If we put $s=\sin^2x/2$ here, then we will obtain that
\[
\pi x^2a(x)=2(1-h)s+2hs(5-20s+16s^2)=2s(1+4h-20hs+16hs^2).
\]
For each $h>0$, the quadratic polynomial in $s$ in the right hand side has the minimum at $s=5/8$, and the minimum value equals $1-9h/4$. We choose $h\leq 4/9$. This implies that $a(x)$ is a density (i.e., $a(x)\geq 0,~ \int a(x)dx=\widehat a(0)=1$).

After the Fourier transform, operator $\mathcal L$ acts as multiplication by $\chi[\widehat a(|k|)-1]$. Thus,  $\mathcal L$ has a continuous spectrum $[-\chi,0]$ and two embedded eigenvalues $\lambda_0=\chi(h-1),~\lambda_1=-\chi$. Each $L_2$-function supported on $[1,2]$ ($[3,\infty)$) is a Fourier representation of an eigenfunction of $\mathcal L$ with the eigenvalue $\lambda_0$ ($\lambda_1$, respectively.)

In order to obtain an operator with infinitely many embedded eigenvalues, consider the function (we will show below that it is a density)
\begin{equation}\label{axx}
a(x)=\frac{1}{\pi x^2}\sum_{j=0}^\infty h_j(\cos2jx-\cos(2j+1)x), \quad \sum h_j=1, \quad x\in R^1\backslash\{0\},
\end{equation}
where $h_j>0$ are chosen in such a way that
\begin{equation}\label{axx1}
h_0>\sum_{j=1}^\infty h_je^{3j},
\end{equation}
i.e., $h_0$ is large enough and $h_j$ decay fast enough. From (\ref{axx1}) it follows that $a(x)$ can be extended by continuity at $x=0$ and that series (\ref{axx}) with complex $x$ converges uniformly on any bounded region of the complex plane. Hence, $a(x)$ is an entire function.

The Fourier image $\widehat{a}=\widehat{a}(|k|)$ of $a(x)$ resembles the function in Fig 1. It is a continuous piece-wise linear function with constant values $c_m=\sum_{j=m}^\infty h_j$ on the intervals $[2m-1,2m],~ m\geq 1, $ and the straight line segments with the slopes $-h_m$ on the intervals $[2m,2m+1],~ m\geq 0. $ The heights of the latter segments are chosen by the continuity condition. Note that $\widehat{a}(0)=c_0=\sum h_j=1$.

Let us rewrite (\ref{axx}) in the form
\[
a(x)=\frac{1}{\pi x^2}\sum_{j=0}^\infty h_j\sin\frac{x}{2}\sin\frac{(4j+1)x}{2}=\frac{1}{\pi x^2}\sin^2\frac{x}{2}[h_0+\sum_{j=1}^\infty h_j\frac{\sin(4j+1)x/2}{\sin x/2}].
\]
The ratio of the sine functions on the right does not exceed $4j+1$. Thus (\ref{axx1}) implies  that $a(x)\geq 0 $. Since $\int a(x)dx=\widehat{a}(0)=1$, function $a(x)$ is a density.

Since the Fourier image $\widehat{a}=\widehat{a}(|k|)$ of $a(x)$ has a constant value $c_m$ on the interval $[2m-1,2m]$, it follows that $\lambda_m=\chi c_m,~m\geq 1,$ are eigenvalues of $\mathcal L$ of infinite multiplicity. The corresponding eigenspaces contain all $L_2$-functions whose Fourier transforms are supported on $[2m-1,2m]$.

The constructed density $a(x)$ satisfies all the requirements of the theorem except the smoothness of its Fourier transform $\widehat{a}(k)$. In order to satisfy  the latter requirement, we replace  $a(x)$ by  $a^{(1)}(x)=c\beta^2(x)a(x)$ where the Fourier transform $\widehat{\beta}(k)$ of $\beta$ is a non-negative function with the support on $|k|\leq\delta<1/5$. We choose $c$ from the condition $\int a^{(1)}(x)dx=1$. Then, obviously, $a^{(1)}(x)$ is a density, is an entire function, and its Fourier transform $\widehat{a^{(1)}}(k)=c\widehat{\beta}(k)*\widehat{\beta}(k)*\widehat{a}(k)$ is infinitely smooth. Since the support of $c\widehat{\beta}(k)*\widehat{\beta}(k)$ belongs to the interval $|k|\leq 2\delta<2/5$, function $\widehat{a^{(1)}}(k)$ has constant values $b_m=bc_m$ on intervals $l_m=[2m-1+2\delta,2m-2\delta], ~m\geq 1.$ Here $|l_m|\geq 1/5$ and
\[
b=c\int_R  \widehat{\beta}(k)*\widehat{\beta}(k)dk=c({\int_R \widehat\beta}(k))^2>0.
\]
Hence $\lambda_m=\chi bc_m,~m\geq 1,$ are eigenvalues of $\mathcal L$ of infinite multiplicity, and the corresponding eigenspaces contain all the $L_2$-functions whose Fourier transforms are supported on $l_m$.

The transition from $d=1$ to arbitrary $d>1$ is simple. If $x=(x_1,...,x_d)$, then one can take
$a(x)=\prod_ia^{(1)}(x_i)$.

\qed

 \section {On the random walk with generator $\mathcal L$.}

The transition density $p(t,x)=P_x\{x(t)=x;~x(0)=0\}$ of the random walk with generator $\mathcal L$
satisfies the equation
\[
p_t=\mathcal Lp,\quad p(0,x)=\delta(x).
\]

We will assume (without loss of generality) that the coefficient $\chi$ in the definition of $\mathcal L$ is equal to one. One can reduce the problem to this case by simple rescaling $t\to \chi t,v(x)\to \chi v(x)$.

Using the Fourier transform, we obtain that
\begin{equation}\label{trd}
p(t,x)=\frac{1}{(2\pi)^d}\int_{R^d}e^{i(k,x)+t(\widehat{a}(k)-1)}dk.
\end{equation}

Condition (\ref{aone1}) implies that the closure of the range range of the function $\widehat{a}(k)-1$ is a segment $[-\alpha, 0]$, where $1\leq\alpha\leq 2$. When $\lambda\in \mathbb C\backslash[-\alpha, 0]$, the Green function $G_\lambda(x-y) $ (the integral kernel of the negative resolvent $G_\lambda:=-(\mathcal L-\lambda)^{-1}$) has the form
  \begin{equation}\label{gg}
 G_\lambda(x-y)=\frac{1}{(2\pi)^d}\int_{R^d}\frac{e^{i(k,x-y)}}{1+\lambda-\widehat{a}(k)}dk.
 \end{equation}

Integrals in both formulas above must be understood as the inverse Fourier transforms in the sense of distributions, but they also can be reduced to convergent integrals.  In particular, the integrand in (\ref{gg}) does not decay at infinity, but we can rewrite the symbol of the operator $G_\lambda$  as
\[
\frac{1}{1+\lambda-\widehat{a}(k)}=\frac{\widehat{a}(k)}{(1+\lambda-\widehat{a}(k))(1+\lambda)}+\frac{1}{1+\lambda}.
\]
We denote the operator with the symbol $\frac{\widehat{a}(k)}{(1+\lambda-\widehat{a}(k))(1+\lambda)}$ by $T_\lambda$ (it acts as multiplication by the symbol in the Fourier images). Then $G_\lambda=\frac{1}{1+\lambda}I+T_\lambda$, where the integral kernel of operator $T_\lambda$ has the form
  \begin{equation}\label{tl}
 T_\lambda(x-y)=\frac{1}{(2\pi)^d}\int_{R^d}\frac{\widehat{a}(k)e^{i(k,x-y)}}{(1+\lambda-\widehat{a}(k))(1+\lambda)}dk, \quad  \lambda\notin [-\alpha,0].
 \end{equation}
 The latter integral converges uniformly in $x,y$ due to (\ref{cond1}), (\ref{cond}). Similarly,

 \begin{equation}\label{pmd}
p(t,x)=e^{-t}\delta(x)+\frac{e^{-t}}{(2\pi)^d}\int_{R^d}e^{i(k,x)}[e^{t\widehat{a}(k)}-1]dk,
\end{equation}
where the integral converges.

 The random process $x(t)$ with the generator $\mathcal L$ is called {\it transient} if
 \begin{equation}\label{trr}
 \lim_{\lambda\to +0} \int_{|y|<r} G_\lambda(x-y)dy  ~~{\rm ~exists ~for ~all}~~ x\in R^d, ~r>0,
 \end{equation}
 and it is called {\it recurrent} in the opposite case. Obviously, the kernel of the operator $G_\lambda=\frac{1}{1+\lambda}I+T_\lambda$ in (\ref{trr}) can be replaced by the kernel of $T_\lambda$. After that, one can easily (using (\ref{cond1}), (\ref{cond})) show that the process $x(t)$ is transient if and only if
 \begin{equation}\label{conv}
\int_{R^d}|\frac{\widehat{a}(k)}{\widehat{a}(k)-1}|dk<\infty.
\end{equation}
From here it also follows that
\begin{equation}\label{rrrr}
 \int_{|y|<r} G_\lambda(x-y)dy\to\infty, ~r>0,\quad  {\rm as} ~~\lambda\downarrow0,
 \end{equation}
uniformly on each compact in $x\in R^d$ if the process is recurrent. Note that in both cases the integrals (\ref{trr}), (\ref{rrrr}) over the whole space aproach infinity as $\lambda\to 0$ since
\[
\int_{R^d} G_\lambda(x-y)dy= \int_0^\infty \int_{R^d}e^{-\lambda t}p(t,x-y)dydt=\frac{1}{\lambda}.
\]

It was shown in the proof of Theorem \ref{ttt} that $\widehat{a}(k)<1$ when $k\neq 0$. Thus the validity of (\ref{conv}) depends only on behavior of $\widehat{a}(k)-1$ near the origin. It follows from (\ref{detb}) that the processes $x(t)$ in the cases of light and moderate tails are recurrent in dimensions $d=1,2$ and are transient if $d\geq 3$. If a slightly stronger version of expansion (\ref{case1}) is valid in the case of heavy tails,  then (see \cite{amv})
$\widehat{a}(k)-1\sim |k|^\gamma$ as $k\to 0$, and therefore the processes with heavy tails are recurrent when $d=1, \gamma\geq 1$, and are
 transient in all the other cases.

The following lemma will be needed for the spectral analysis of operator $H$.
\begin{lemma} \label{l1}
Functions $p(t,x)-e^{-t}\delta(x)$ for $t>0$, $T_\lambda(x)$ for $\lambda>0$, and $T_0(x)$ when the process $x(t) $ is transient are strictly positive.
\end{lemma}
 {\bf Proof.}  From (\ref{cond}) it follows that function $a(x)$ is uniformly continuous and  bounded: $|a(x)|\leq A:=\frac{1}{(2\pi)^d}\int_{R^d}|\widehat{a}(k)|dk$. Put
\begin{equation}\label{an0}
 a_n(x)=(a*a*...*a)(x),
 \end{equation}
 where there are exactly $n$ convolution factors on the right. Using the last relation in (\ref{aone1}), one can easily justify by induction that all the functions $a_n(x)$ are continuous and bounded by the same constant $A$. The latter is useful to justify the convergence of the series (\ref{pb0}) below.

 Due to the symmetry of $a(x)$, we have
 \[
 a_2(0)=\int_{R^d}a(y)a(-y)dy=\int_{R^d}a^2(y)dy>0.
 \]
 Thus there exists $\delta>0$ such that $a_2(x)>0$ when $|x|<\delta$. Since $a(x)\geq 0$, from the definition of the convolution it follows that $a_n(x)\geq 0$ for all $n, x$, and  $a_{2n}(x)>0$ when $|x|<n\delta$. The latter two facts together with (\ref{pmd}) immediately imply that
  \begin{equation}\label{pb0}
 p(t,x)-e^{-t}\delta(x)=e^{-t}\sum_{n=1}^\infty \frac{t^n}{n !}a_n(x)>0, \quad t>0.
\end{equation}

The positivity of the first function in the statement of the lemma is proved. The positivity of the second one could be proved similarly. Another option is to note that
  \begin{equation}\label{ltl}
T_\lambda(x)=\int_0^\infty [p(t,x)-e^{-t}\delta(x)]e^{-\lambda t}dt, \quad \lambda>0.
\end{equation}
The formula above follows from (\ref{tl}), (\ref{pmd}) and the uniform boundedness of $p(t,x)-e^{-t}\delta(x)$, which is one of the consequences of (\ref{pb0}). The positivity of the integrand in (\ref{ltl}) implies that $T_\lambda(x)>0$ when $\lambda>0$. Moreover, the right-hand side in (\ref{ltl}) increases monotonically when $\lambda\to +0$. Hence, if the limiting function $T_0(x)$ exists, it is also positive.

\qed

 \section {Spectral analysis of operator $H$}
We will continue to assume that $\chi=1$.

Operator $H=\mathcal L+v(x)$ is a sum of a convolution operator and an operator of multiplication by a compactly supported potential. There is a duality between these two terms. If $H$ is rewritten in the form $H=(\mathcal L+1)+(v(x)-1)$, then after the Fourier transform, the first term becames an operator of multiplication by function $\widehat{a}(k)$ decaying at infinity (see (\ref{cond1})) and the second term becomes a convolution operator. If $v$ is smooth, then the spectrum of the second term is given by Lemma \ref{laaa}, and it coincides with the segment $[-1, \max v(x)-1]$. These arguments lead to the following statement.
\begin{theorem}
Let the potential $v$ be continuous. Then the essential spectrum of operator $H$ contains segments $[-a,0]$ and $[-1, \max v(x)-1]$.
\end{theorem}

{\bf Proof.} Let $\lambda_0\in [-1, \max v(x)-1]$ and $ v(x_0)-1=\lambda_0$. Let
$$
\psi_\varepsilon(x)=(\pi\varepsilon)^{-d/4}e^{-\frac{|x-x_0|^2}{2\varepsilon}}.
$$
Then $\|\psi_\varepsilon(x)\|_{L^2}=1$ and $\|( v(x)-1)\psi_\varepsilon(x)\|_{L^2}\to 0$ as $\varepsilon\to 0$. Moreover,
$$
|\widehat{\psi}_\varepsilon(k)|=(4\pi\varepsilon)^{d/4}e^{-\frac{\varepsilon|k|^2}{2}}.
$$
Hence $\|\widehat{a}\widehat{\psi}_\varepsilon\|_{L^2}\to 0$ as $\varepsilon\to 0$ due to the decay of $\widehat{a}(k)$ at infinity, i.e.,
$$
\|( \mathcal L+1)\psi_\varepsilon(x)\|_{L^2}\to 0 ~~~ {\rm as}~~~ \varepsilon\to 0.
$$
 Thus $\psi_\varepsilon(x)$ is a Weyl sequence, and therefore $\lambda_0$ belongs to the essential spectrum of $H$. The segment $[-a,0]$ is treated similarly due to the duality discussed above.

 \qed

We would like to consider operators where the potential can produce only discrete spectrum on positive semiaxis. Thus we will assume that $v(x)\leq 1$, and often that $v(x)< 1$. As it was mentioned in the introduction, this restriction appears naturally in contact models. On the other hand, we would like to introduce a sequence of potentials $v_R(x)$ whose action becomes stronger as $R\to \infty$. Thus we will consider $v_R(x)=v(x/R)$. We do not need for the potential to be smooth. It can be only measurable (and bounded), but in some cases we assume that it has a positive limiting value at the origin or is continuous.

\begin{theorem}\label{sss}
Let the potential $v$ be continuous, $0\leq v(x)\leq 1-\delta$ with some $\delta>0$, and $\int_{R^d}v(x)dx>0$. Then the following statements are valid.

(1) The spectrum of the operator $H=\mathcal L+v(x)$ on the positive semi-axis consists of at most a countable set of eigenvalues of finite multiplicity with the only possible limiting point at $\lambda=0$. These eigenvalues are located on the interval $(0,1-\delta].$

(2) If the underlying process is recurrent, then operator $H=\mathcal L+v(x)$ has at least one positive eigenvalue.

(3) Let
\begin{equation}\label{hr}
H=H(R)=\mathcal L+v(x/R),
\end{equation}
and let underlying process be transient. Then $H$ does not have positive eigenvalues when $R$ is small enough. If, additionally, $v(0)>0$, then positive eigenvalues exist for large enough $R$.


(4) In all the cases, the largest positive eigenvalue $\lambda_0(R)$ of operator (\ref{hr}) (if positive eigenvalues exist) is simple and the corresponding eigenfunction is positive. Function $\lambda_0(R)$ is a monotone function of $R$ if $v(x)$ depends monotonically on $|x|$.
 \end{theorem}

 {\bf Proof.}  In all the cases, we will assume that the operator $H$ has form (\ref{hr}). One can choose $R=1$ in the proof of the first two statements.

 Formula $\psi =\sqrt {v_R} u$ establishes a one-to-one correspondence between the eigenfunctions $u$ of $H$ with positive eigenvalues $\lambda=\lambda_i>0$ and the solutions of the problem
 \begin{equation}\label{gt}
 (I- \sqrt {v_R}G_\lambda\sqrt {v_R})\psi=0, \quad  \psi\in L^2_{{\rm com}}, \quad \lambda>0.
\end{equation}
We put here $G_\lambda=\frac{1}{1+\lambda}I+T_\lambda$ and arrive at
\[
((1-\frac{ v_R}{1+\lambda})I-\sqrt {v_R}T_\lambda\sqrt {v_R})\psi=0, \quad  \psi\in L^2_{{\rm com}},
\]
or, equivalently,
 \begin{equation}\label{tt}
 (I- \sqrt {w}T_\lambda\sqrt {w})\psi_1=0, \quad  w=\frac{(1+\lambda)v_R}{1+\lambda- v_R}\geq 0,   \quad \psi_1=\sqrt {w}u\in L^2(D),
\end{equation}
where $D$ is the support of $w(x)$.

We will consider $\sqrt {w}T_\lambda\sqrt {w}$ as an operator in $L^2(D)$. Since the symbol $\frac{\widehat{a}(k)}{(1+\lambda-\widehat{a}(k))(1+\lambda)}$  of $T_\lambda$ is integrable, its inverse Fourier transform is bounded and the integral kernel of the operator $\sqrt {w}T_\lambda\sqrt {w}$ is bounded and compactly supported. Thus operator $\sqrt {w}T_\lambda\sqrt {w}$ is compact (this is the reason to consider (\ref{tt}) instead of (\ref{gt})). Obviously, $\sqrt {w}T_\lambda\sqrt {w}$
depends analytically on $\lambda$ when $\lambda>0$, and
\begin{equation}\label{nnn}
\|T_\lambda\|_{L^2}\to 0 ~~~{\rm as}~~~\lambda\to\infty.
\end{equation}
Thus the operator in the left-hand side of (\ref{tt}) is invertible when $\lambda\gg 1$. The analytic Fredholm theorem implies the validity of the first statement of the theorem except the assertion that the eigenvalues belong to $(0,1-\delta].$
The latter follows from the fact that $\mathcal L\leq 0$ (and therefore $H\leq \sup v(x)$).

For a more detailed analysis of the set $\{\lambda_i\}$, consider the eigenvalues $\{\mu_j(\lambda,R)\}, ~\lambda>0, R>0, $ of the operator $\sqrt {w}T_\lambda\sqrt {w}$. Since $T_\lambda(x-y)$ is strictly positive (due to Lemma \ref{l1}), the Perron-Frobenius theorem implies that, for each $\lambda, R>0$, the operator $\sqrt {w}T_\lambda\sqrt {w}$ has the largest positive simple eigenvalue $\mu_0$ (the ground energy) with a positive eigenfunction (the ground state), and $|\mu_j|<\mu_0, ~j>0$. Obviously, $ \mu_j(\lambda,R)\to +0$ as $j\to\infty$. Equation  (\ref{tt}) implies that the eigenvalues $\lambda_i>0$ of $H$ are defined by the equations
\begin{equation}\label{muj}
\mu_j(\lambda,R)=1, j\geq 0.
\end{equation}

From (\ref{pb0}) and (\ref{ltl}) it follows that
  \begin{equation*}
T_\lambda(x)=\sum_{n=1}^\infty \frac{a_n(x)}{(1+\lambda)^{n+1}}
\end{equation*}
is a strictly decaying function of $\lambda$ when $\lambda>0$. From (\ref{tt}) it follows that $w=1+\frac{v_R}{1+\lambda- v_R}$ is a strictly decaying function of $\lambda$ when $v_R\neq 0$ and that $w\to v(0)$ monotonically as $R\to\infty$ if $v(x)$ is monotone in $|x|$. The positivity of the integral kernel of operator  $\sqrt {w}T_\lambda\sqrt {w}$ and the positivity of the ground state together with the monotonicity properties discussed above imply that the ground state
\[
\mu_0(\lambda, R)=\min_{\psi:\|\psi\|=1}(\sqrt {w}T_\lambda\sqrt {w}\psi,\psi)
\]
is a strictly decaying function of $\lambda,~\lambda>0$, and that it is increases with $R$ if  $v(x)$ is monotone in $|x|$.

Recall that $\mu_0(\lambda, R)\to 0$ as $\lambda\to \infty$ (due to (\ref{nnn})). The monotonicity of $\mu_0(\lambda, R)$ in $\lambda$ implies the existence of the limit $\mu_0(+0, R)=\lim_{\lambda\to +0}\mu_0(\lambda, R)$. From (\ref{muj}) it follows that $H$ has positive eigenvalues if and only if $\mu_0(+0, R)> 1$, and the largest eigenvalue $\lambda_0(R)>0$ of $H$ is defined by the equation
\begin{equation}\label{mu0}
\mu_0(\lambda,R)=1.
\end{equation}

Hence from the Perron-Frobeneus theorem and the monotonicity of $\mu_0$ it follows that the eigenvalue $\lambda_0(R)>0$ is simple and monotone in $R$ when $v(x)$ is monotone in $|x|$. Moreover, the eigenfunction $u$ of $H$ with the eigenvalue $\lambda_0(R)$ satisfies the equation $(\mathcal L-\lambda_0+v_R)u=0$, and $\psi_1=\sqrt wu$
is the the ground state for $\sqrt wT_\lambda\sqrt w$. Thus $u=G_\lambda\psi_1>0$ since $\psi_1\geq0$ and since the integral kernel of operator $G_\lambda=\frac{1}{1+\lambda}+T_\lambda$ is strictly positive due to Lemma \ref{l1}.
 This completes the proof of the last statement of the theorem. It remains to prove statements two and three. Prior to starting the proofs, let us note that the Perron-Frobeneus theorem implies that $\mu_0(\lambda, R)=\|\sqrt wT_\lambda\sqrt w\|$. In particular, $H$ has positive eigenvalues if and only if
 \begin{equation}\label{xxx}
 \lim_{\lambda\to +0}\|\sqrt wT_\lambda\sqrt w\|>1.
 \end{equation}

Assume that the process $x(t)$ is recurrent. Then (\ref{rrrr}) implies that $\|\sqrt {w}T_\lambda\sqrt {w}\|\to\infty$ as $\lambda\to0$. This proves the second statement of the theorem.

Let the underlying random process be transient, i.e., (\ref{conv}) holds. The latter relation implies that the integral kernel $T_\lambda(x-y)$ of the operator $T_\lambda$ is continuous and converges uniformly in $x,y$ to a continuous function $T_0(x-y)$ as $\lambda\to +0$.
Then the operators $\sqrt {w}T_\lambda\sqrt {w}$ converge in the operator norm to the integral operator $Q$ with the kernel
\[
\sqrt{\frac{v_R(x)}{1- v_R(x)}}T_0(x-y)\sqrt{\frac{v_R(y)}{1- v_R(y)}}, \quad T_0(x-y)=\frac{1}{(2\pi)^d}\int_{R^d}\frac{\widehat{a}(k)e^{ik(x-y)}}{1-\widehat{a}(k)}dk\in C.
\]
Due to (\ref{xxx}), $H$ has positive eigenvalues if and only if $\|Q\|>1$.

Since function $|T_0(x)|$ is bounded due to (\ref{conv}), the integral kernel $Q(x,y)$ of operator $Q$ has the estimate $|Q(x,y)|<\sqrt{\frac{v_R(x)}{1- v_R(x)}}\sqrt{\frac{v_R(y)}{1- v_R(y)}}$. Hence
\[
\|Q\|^2\leq \int |Q(x,y)|^2dxdy\leq C(\int \frac{v(x/R)}{1- v(x/R)}dx)^2\to 0 \quad {\rm as } ~~R \to 0.
\]
Thus $H$ does not have positive eigenvalues when $R$ is small enough.

In order to prove that the eigenvalues exist when $R$ is large enough, we need to show that $\|Q\|>1$ for $R\gg 1$. Consider $\psi_\varepsilon(x)=\varepsilon^{d/4}\pi^{-d/4}e^{-\varepsilon|x|^2/2}$. Then $\|\psi_\varepsilon\|=1 $ and $\widehat{\psi}_\varepsilon(k)=(\varepsilon\pi)^{-d/4}e^{-|x|^2/2\varepsilon}$. Obviously, $\widehat{\psi}_\varepsilon^2\to\delta(k)$ as $\varepsilon\to+0$ and
\[
(T_0\psi_\varepsilon,\psi_\varepsilon)=(2\pi)^{-d}\int_{R^d}\frac{\widehat{a}(k)}{1-\widehat{a}(k)}\widehat{\psi}_\varepsilon^2(k)dk\to\infty  \quad {\rm  as} ~~ \varepsilon\to+0
\]
since $\widehat{a}(0)=1$. We assume that $v(0)>0$ and choose $\varepsilon=\varepsilon_0>0$ such that the left-hand side above exceeds $2v^{-2}(0)$ when $0<\varepsilon\leq\varepsilon_0$. We have
\[
(Q\psi_{\varepsilon_0},\psi_{\varepsilon_0})=(2\pi)^{-d}\int_{R^d}\frac{\widehat{a}(k)}{1-\widehat{a}(k)}|\widehat{v_R\psi}_{\varepsilon_0}|^2(k)dk,
\]
where $\widehat{v_R\psi}_{\varepsilon_0}$ is the Fourier transform of $v(x/R)\psi_{\varepsilon_0}(x)$. The latter product converges to $v(0)\psi_{\varepsilon_0}(x)$ in $L_1$ as $R\to\infty$, and therefore $\widehat{v_R\psi_{\varepsilon_0}}\to v(0)\widehat{\psi_{\varepsilon_0}}$ uniformly as $R\to\infty$. Hence (\ref{conv}) implies that
\[
(Q\psi_{\varepsilon_0},\psi_{\varepsilon_0})\to v^2(0)(T_0\psi_{\varepsilon_0},\psi_{\varepsilon_0})>2 \quad {\rm as} ~~R\to\infty.
\]
This proves the existence of the eigenvalues when $R\gg 1$.

\qed

\section {Asymptotics of the Green function at infinity}
We will obtain here the asymptotics of the Green function $G_\lambda(x),~\lambda>0,~|x|\to \infty,$ for processes with ultra light tails, i.e., under the assumption that
 \begin{equation}\label{ult}
a(x)\leq Ce^{-|x|^\alpha},~~ \alpha>1.
\end{equation}
We will also assume that the same estimate is valid for the gradient of $a(x)$:
\begin{equation}\label{ulta}
|\nabla a(x)|\leq Ce^{-|x|^\alpha}.
\end{equation}

From (\ref{gg}), (\ref{tl}) it follows that $G_\lambda(x),~\lambda>0,$ can be rewritten in the form (compare to (\ref{pb0}))
 \begin{equation}\label{gg1}
G_\lambda(x)=\frac{\delta(x)}{1+\lambda}+\sum_{n=1}^\infty\frac{a_n(x)}{(1+\lambda)^{n+1}},
\end{equation}
where $a_n$ is the convolution of $n$ copies of $a(x)$, see (\ref{an0}). The asymptotics of $G_\lambda$ at infinity will be expressed in terms of the moment generating function
\begin{equation}\label{mgf}
Ee^{(\nu,Y)}:=\int_{R^d}e^{(\nu,y)}a(y)dy,
\end{equation}
where $Y$ is a random variable with the density $a(\cdot)$. From (\ref{ult}) it follows that the function (\ref{mgf}) is analytic in $\nu$ and positive when $\nu\in R^d$. Thus it can be rewritten as
\begin{equation}\label{hess11}
\int_{R^d}e^{(\nu,y)}a(y)dy=e^{H(\nu)}, \quad \nu\in R^d,
\end{equation}
where ${H(\nu)}, ~ \nu\in R^d,$ is analytic and real-valued.
\begin{lemma}\label{hes}
Function $H(\nu)$ has the following properties: $H(0)=0,~ H(\nu)$ is even ($H(\nu)=H(-\nu)$), and it is strictly convex with
\begin{equation}\label{hess1}
B(\nu):={\rm Hess}H(\nu)=[\frac{\partial^2H(\nu)}{\partial\nu_i\partial\nu_j}]>0.
\end{equation}
\end{lemma}

{\bf Proof.} Relation (\ref{aone1}) implies that $H(0)=0$. The symmetry of $H$ is a consequence of the symmetry of $a$. Let us show (\ref{hess1}). Consider the density
\begin{equation}\label{anu}
a_\nu(y)=e^{-H(\nu)}a(y)e^{(\nu,y)}, \quad \nu\in R^d.
\end{equation}
Obviously, $a_\nu(y)\geq 0$ and $\int_{R^d}a_\nu(y)dy=1$, i.e., $a_\nu(y)$ is a density. Since
\[
\int_{R^d}a_\nu(y)e^{(k,y)}dy=e^{H(\nu+k)-H(\nu)},
\]
we have
\begin{equation}\label{expnu}
\int_{R^d}ya_\nu(y)dy=\nabla_k[e^{H(\nu+k)-H(\nu)}]_{k=0}=\nabla H(\nu).
\end{equation}
Similarly,
\[
\int_{R^d}[y_iy_j]a_\nu(y)dy=[\frac{\partial H(\nu)}{\partial\nu_i}\frac{\partial H(\nu)}{\partial\nu_j}]+[\frac{\partial^2H(\nu)}{\partial\nu_i\partial\nu_j}].
\]
The last two relations imply that $B(\nu)$ is the covariance of the process with the density $a_\nu(y)$, and therefore $B(\nu)>0$.

\qed

Denote by $H^*(p)$ the Legendre transform of the function $H(\nu)$, i.e.,
\begin{equation}\label{hessh}
H^*(p)=\max_\nu[(p,\nu) -H(\nu)],
\end{equation}
where the maximum can be equal to infinity. If a finite maximum exists, then
\begin{equation}\label{hess2}
H^*(p)=(p,\nu^*(p)) -H(\nu^*(p)),
\end{equation}
where $\nu^*(p)$ is the solution of the equation
\begin{equation}\label{hess2a}
\nabla H(\nu)=p\in R^d.
\end{equation}
Let us describe the set $\{p\}$ where the critical point $\nu^*(p)$ exists and function $H^*(p)$ is finite.

Denote by $S_0$ the set of points $y\in R^d$ where $a(y)>0$.  Recall that  $Y$ is symmetric (see (\ref{aone1})), and therefore the set $S_0$ is symmetric. Let $\mathcal F$ be the convex hull of the set $S_0$, i.e.,
\[
\mathcal F=\bigcap_{\dot{\nu}}L(\dot{\nu}), \quad \dot{\nu}=\frac{\nu}{|\nu|}\in S^{d-1},
\]
where $L(\dot{\nu})=\{y:|(y,\dot{\nu})|<s^+(\dot{\nu}) \}$ is the minimal layer containing $S_0$ and orthogonal to $\dot{\nu}$. Obviously, vectors $\nu$ for which $s^+(\dot{\nu})=\infty$ combined with the zero vector form a linear subspace $M\subset R^d$. The set $\mathcal F$ is a cylinder: it is translation invariant with respect to vectors from $M$, and the cross-section of $\mathcal F$ orthogonal to $M$ is a convex bounded set.
\begin{lemma}\label{hes3} The following statements are valid:

1) For each $p\in\mathcal F$, equation (\ref{hess2a}) has a unique solution $\nu=\nu^*(p)$, and function  (\ref{hessh}) is finite. Its values are given by (\ref{hess2}).

2) $H^*(p)=\infty$ if $p\notin\mathcal F$.

3) $H^*(p)\to\infty$ as $p$ approaches the boundary of $\mathcal F$ or goes to infinity. Moreover,  $\frac{H^*(p)}{|p|}\to\infty$ as $|p|\to\infty$.
\end{lemma}
{\bf Remark.} This lemma is not valid in the lattice case.

{\bf Proof.} We rewrite (\ref{hess11}) in the form
\begin{equation}\label{d-1}
H(\nu)=\ln\int_{R^d}e^{|\nu|(\dot{\nu},y)}a(y)dy=
\ln\int_{-s^+}^{s^+}e^{|\nu|\tau}a_{\dot{\nu}}(\tau)d\tau, \quad {\rm where} \quad a_{\dot{\nu}}(\tau)=
\int_{(\dot{\nu},y)=\tau}a(y)dy',
\end{equation}
and where $dy'$ is the  volume element in $R^{d-1}$.
Here $s^+=s^+(\dot{\nu})$, $a_{\dot{\nu}}(\tau)=0$ when $|\tau|>s^+$,  and $s^+$ is a limiting point of the set $\{\tau: a_{\dot{\nu}}(\tau) >0 \}$.
  The latter property of $a_{\dot{\nu}}$ and (\ref{d-1}) imply that $H(\nu)>s'|\nu|$ for each $s'<s^+$ if $|\nu|$  is large enough. On the other hand, $|(p,\dot{\nu})|<s^+$ for $p\in \mathcal F$. Thus
\begin{equation}\label{inff}
(p,\nu)-H(\nu)\to -\infty
\end{equation}
 when $|\nu|\to\infty$ and $p\in \mathcal F$ is fixed. Hence, function (\ref{hessh}) is finite. Relations (\ref{inff}) and (\ref{hess1}) imply that equation (\ref{hess2a}) has a unique solution $\nu=\nu^*(p)$ when $p\in \mathcal F$ and that (\ref{hess2}) is valid.

 If $p\notin \mathcal F$, then there is a unit vector $\dot{\nu_0}$ such that $s^+=s^+(\dot{\nu_0})<\infty$ and $(p,\dot{\nu_0})\geq s^+$. On the other hand, (\ref{hess1}) implies that $H(\sigma\dot{\nu_0})<\ln\frac{Ce^{s^+\sigma}}{\sigma}$, where $C$ is an upper bound for $|a_{\dot{\nu}}(\tau)|$. Thus
 \[
 (p,\sigma\dot{\nu_0})-H(\sigma\dot{\nu_0})\geq s^+\sigma-\ln\frac{Ce^{s^+\sigma}}{\sigma}=\ln\sigma-\ln C\to\infty \quad {\rm as} \quad \sigma\to\infty,
 \]
 and therefore the right-hand side in (\ref{hessh}) is infinite. It remains to prove the last statement of Lemma \ref{hes3}.

Let  $p\to p'\in \partial\mathcal F$. As above, there is $\dot{\nu_0}$ such that $s^+=s^+(\dot{\nu_0})<\infty$ and $(p',\dot{\nu_0})= s^+$. Then
\begin{eqnarray*}
H^*(p)\geq \max_{\sigma\geq 1}[(p,\sigma\dot{\nu_0})-H(\sigma\dot{\nu_0})]=\max_{\sigma\geq 1}[(p-p',\sigma\dot{\nu_0})+s^+\sigma -H(\sigma\dot{\nu_0})]\\  \geq \max_{\sigma\geq 1}[-\varepsilon\sigma+\ln\sigma-\ln C],
\end{eqnarray*}
where $\varepsilon=\varepsilon(p)=(p-p_0,\dot{\nu_0})$. The last maximum goes to infinity as $\varepsilon\to 0$. Thus, $H^*(p)\to\infty$ as $p\to p_0$.

Let $|p|\to \infty$ and $\nu'=np/|p|, n>1$. From (\ref{hessh}) it follows that
\[
H^*(p)\geq (p,\nu')-H(\nu')\geq n|p|-C(n), \quad {\rm where} \quad C(n)=\max_{|\nu|\leq n}H(n).
\]
Since $n$ is arbitrary, the latter estimate implies that $\frac{H^*(p)}{|p|}\to \infty$ as $|p|\to \infty$.

\qed

Consider the
``phase" function
\[
S=S(\tau,\theta,\lambda)=\tau[ H^*(\frac{\theta}{\tau})+\ln(1+\lambda)], \quad
0<\tau<\infty,
\]
where $\theta=\frac{x}{|x|},~~\varepsilon\leq \lambda\leq 1/\varepsilon$.
Due to Lemma \ref{hes3}, this function is smooth in $(\tau, \theta, \lambda)$ when $\frac{1}{\tau}<s^+(\theta)$ and is equal to infinity if $\frac{1}{\tau}\leq s^+(\theta)<\infty$. Moreover, in both cases, $S\to\infty$ as $\frac{1}{\tau}\to s^+(\theta)$. Note also that $S\to\infty$ as $\tau\to \infty$. Thus, for each $\theta$ and $\lambda$, function
$S$ achieves its minimum value. This value is positive since $S$ is strictly positive. Note also that direct calculations imply that
 \[
 S_{\tau\tau}=\frac{1}{\tau}<{\rm Hess}H^*(p)p,p>, \quad  p=\frac{\theta}{\tau}, \quad |p|<s^+(\theta),
 \]
and therefore $S_{\tau\tau}> 0$. Thus for each $\theta$ and $\lambda$, there is a unique positive point $\tau=\tau^0(\theta,\lambda)$ where $S$ has the absolute minimum value, and this point depends smoothly on $\theta$ and $\lambda$. Let
\begin{equation}\label{fi}
\varphi=\varphi(\theta,\lambda):=S|_{\tau=\tau^0(\theta,\lambda)}>0.
\end{equation}
\begin{theorem}\label{tgr}
There exists a smooth function $f=f(\theta, \lambda)>0$ such that for each $\lambda'>0,$
\begin{equation}\label{asg}
G_\lambda(x)= f(\theta, \lambda)|x|^{(1-d)/2}e^{-|x|\varphi(\theta,\lambda)}(1+o(1)), \quad 0<\lambda'\leq \lambda\leq 1/\lambda', \quad |x|\to\infty,
\end{equation}
where $\varphi$ is defined in (\ref{fi}).
\end{theorem}
In order to prove the theorem, we will need the following two lemmas.
\begin{lemma}\label{ann}
The following asymptotics holds uniformly in $y$ when $\frac{|y|}{n}\leq b<\infty$:
\begin{equation}\label{abc}
a_n(y)= \frac{e^{-nH^*(\frac{y}{n})}}{(2\pi n)^{d/2}\sqrt{\det B(\nu^*(\frac{y}{n}))}}(1+o(1)),
 \quad n\to\infty.
\end{equation}
\end{lemma}
{\bf Remark.} Logarithmic asymptotics for $a_n(y)$ and local asymptotics for small $\frac{|y|}{n}$ can be found in \cite{pet}, \cite{wfr}.

{\bf Proof.} Consider the sum $S_{n,\nu}=X_{1,\nu}+...+X_{n,\nu}$ of i.i.d.r.v. with the density $a_\nu(y)$ (see (\ref{anu})). Due to the CLT,
\[
\frac{S_{n,\nu}-EX_{1,\nu}}{\sqrt n} ~\begin{array}{c} {\rm law} \\ \longrightarrow \\  \\  \end{array} ~N(0, B(\nu)), \quad n\to\infty.
\]
In fact we need the asymptotics of the density $a_{n,\nu}(y)$ of the sum $S_{n,\nu}$. Due to the local CLT \cite{Fel},
\[
a_{n,\nu}(nEX_{1,\nu})=\frac{1}{(2\pi n)^{d/2}\sqrt{\det B(\nu)}}(1+o(1)), \quad n\to\infty,
\]
uniformly in $\nu$ in every ball $|\nu|\leq \beta<\infty$. The local limit theorem for densities holds under the condition that
\[
|\widehat{a_{\nu}}|^m \in L^1(R^d)
\]
for some $m>0$. This condition holds in our case due to (\ref{ult}), (\ref{ulta}).
We combine the local limit theorem  for $a_{n,\nu}(nEX_{1,\nu})$ with (\ref{expnu}) and obtain
\[
a_{n,\nu}(n\nabla H(\nu))=\frac{1+o(1)}{(2\pi n)^{d/2}\sqrt{\det B(\nu)}}, \quad |\nu|\leq \beta, \quad n\to\infty.
\]

Function $a_{n,\nu}$ is the convolution of $n$ copies of (\ref{anu}), and therefore
\[
a_{n,\nu}(y)=e^{-nH(\nu)}a_n(y)e^{(\nu,y)}.
\]
Thus,
\[
a_{n}(n\nabla H(\nu))=e^{nH(\nu)-(\nu,y)}a_{n,\nu}(n\nabla H(\nu))=\frac{e^{nH(\nu)-(\nu,y)}}{(2\pi n)^{d/2}\sqrt{\det B(\nu)}}(1+o(1)),
\]
where $|\nu|\leq \beta, ~ n\to\infty$. It remains only to choose $\nu=\nu^*(\frac{y}{n})$.

\qed

\begin{lemma}\label{ann2}
There exists a constant $c$ such that
\begin{equation}\label{abc2}
|a_n(y)|\leq e^{n(c-\frac{1}{2}(\frac{|y|}{n})^\alpha)} \quad  when \quad  \frac{|y|}{n}\geq 1.
\end{equation}
\end{lemma}
{\bf Remark.} This estimate shows that $a_n(y)$ decays at infinity somewhat slower, but similarly to $a(y)$ (see (\ref{ult})). The estimate will be used when $\frac{|y|}{n}$ is large enough.

{\bf Proof.} We write $a(y)$ in the form $a(y)=e^{-\frac{1}{2}|y|^\alpha}b(y)$, Then $b\in L^1$ due to (\ref{ult}). Consider the integrand in the convolution of $n$ copies of $a(y)$. We combine together all the exponents and all the factors with function $b$. Then we estimate the product of the exponents by its maximum value. This implies that
\[
|a_n(y)|\leq(\max e^{-\frac{L}{2}})|b_n(y)|, \quad L=|y-y^{(1)}|^\alpha+|y^{(1)}-y^{(2)}|^\alpha+...+|y^{(n-1)}|^\alpha,
\]
where $b_n$ is the convolution of $n $ copies of the function $b$. It is obvious that $L$ has the minimum value when all its terms are equal, i.e., $L\geq n(\frac{|y|}{n})^\alpha$. It remains to note that $|b_n(y)|\leq C^n=e^{nc}$, where $C=\|b\|_{L^1},~ c=\ln C$.

\qed

{\bf Proof of Theorem \ref{tgr}.} We represent $G_\lambda(x),~x\neq0,$ as the sum $G_\lambda(x)=G_1+G_2+G_3,~x\neq0,$ where $G_1$ contains terms from the right-hand side of (\ref{gg1}) with small values of $n$, $G_2$ contains terms with large values of $n$, and $G_3$ contains terms with  intermediate values of $n$. We do not need to worry about the delta-function in the right-hand side of (\ref{gg1}) since we assume that $x\neq 0$. To be more exact, the splitting of $G_\lambda(x)$ depends on the values of $x$ and $\lambda$, and is defined as follows.

Let $s\geq 1$ be an arbitrary number such that
\begin{equation}\label{abc3}
\frac{1}{2}t^\alpha -2\varphi t-c+\ln(1+\lambda)> 0 \quad {\rm for} ~~~ t\geq s,
\end{equation}
where $\varphi=\varphi(\theta,\lambda)$ and $c>0$ are defined in (\ref{fi}) and Lemma \ref{ann2}, respectively. Function $G_1$ is the sum of terms from (\ref{gg1}) with $n$ such that  $\frac{|x|}{n}\geq s$. Function $G_2$ consists of terms with $\frac{|x|}{n}\leq s_1$, where $s_1$ is an arbitrary number from the interval $(0,\frac{\ln(1+\lambda)}{\varphi})$. The inequalities $s_1<\frac{|x|}{n}< s$ hold for the terms in $G_3$. Let us show that $G_1$ and $G_2$ do not contribute to the asymptotics of $G_\lambda(x)$.

Since $s\geq 1$, from Lemma \ref{ann2} and (\ref{abc3}) it follows that
\[
G_1\leq C\sum_{n\leq |x|/s}e^{n[c-\frac{1}{2}(\frac{|x|}{n})^\alpha -\ln(1+\lambda)]}\leq C\sum_{n\leq |x|/s}e^{-n(2\varphi\frac{|x|}{n})}= C\sum_{n\leq |x|/s}e^{-2\varphi |x|}\leq C_1|x|e^{-2\varphi |x|}.
\]
This term is exponentially smaller than the right-hand side in (\ref{asg}).

Lemmas \ref{ann} and \ref{hes} imply that $a_n(x)$ are uniformly bounded when $\frac{|x|}{n}\leq s_1$. Hence
\[
G_2\leq  C\sum_{n\geq |x|/s_1}e^{-n\ln(1+\lambda)}\leq C_1e^{-|x|\ln(1+\lambda)/s_1}\leq C_1e^{-|x|(\varphi+\gamma)},
\]
where $\gamma>0$ since $s_1<\frac{\ln(1+\lambda)}{\varphi}$. Hence $G_2$ is also exponentially smaller than the right-hand side in (\ref{asg}).

Let us evaluate $G_3$ now. We have
\[
G_3=\sum_{n:~s_1<\frac{|x|}{n}< s}\frac{a_n(x)}{(1+\lambda)^{n+1}}=\sum_{n:~s_1<\frac{|x|}{n}< s}n^{-d/2}g(\frac{x}{n})e^{-n[H^*(\frac{x}{n})+\ln(1+\lambda)]}(1+o(1)),
\]
where $g(\frac{x}{n})=C[\det B(\nu^*(\frac{x}{n})]^{-1/2}$. We introduce $\tau_n=\frac{n}{x}$ and rewrite the last formula in the form
\[
G_3=\sum_{n:~1/s<\tau_n< 1/s_1}|x|^{-d/2}h(\tau_{n})e^{-|x|S(\tau_n,\theta,\lambda)}(1+o(1)),
\]
where $h(\tau):=\tau^{-d/2}g(\frac{1}{\tau})$ is infinitely smooth on the segment $[1/s, 1/s_1]$. Since $n>|x|/s$ here, the remainder terms $o(1)$ in the formula above vanish uniformly in $n$ when $|x|\to\infty$ (see Lemma \ref{ann}). Besides, all the terms in the right-hand side are positive. Hence, it is enough to prove the statement of the theorem for
\begin{equation}\label{1221}
u:=\sum_{n:~1/s<\tau_n< 1/s_1}|x|^{-d/2}h(\tau_{n})e^{-|x|S(\tau_n,\theta,\lambda)}.
\end{equation}

Recall that $S$ has the absolute minimum value $\varphi(\theta,\lambda)$ when $\tau=\tau^{(0)}$. We fix an arbitrary $\varepsilon>0$ and choose a $\delta>0$ such that $|S_\tau|<\varepsilon$ when $|\tau-\tau^{(0)}|\leq\delta$. Since $S(\tau,\theta,\lambda)-\varphi(\theta,\lambda)>\gamma_1>0$ when $|\tau-\tau^{(0)}|\geq \delta$ and the number of terms in (\ref{1221}) has order $O(|x|)$, formula (\ref{1221}) can be rewritten as follows
\begin{equation}\label{1222}
u=\sum_{n:~|\tau_n-\tau^{(0)}|<\delta}|x|^{-d/2}h(\tau_{n})e^{-|x|S(\tau_n,\theta,\lambda)}+
O(e^{-|x|[\varphi(\theta,\lambda)+\gamma_1]}).
\end{equation}
In order to obtain the asymptotics of $u$, we will show that function (\ref{1222}) is close to the integral
\[
I=|x|^{(2-d)/2}\int_{|\tau-\tau^{(0)}|<\delta}h(\tau)e^{-|x|S(\tau,\theta,\lambda)}d\tau
\]
when $|x|\to\infty$. Then we will apply the Laplace method to get the asymptotics of $I$ as $|x|\to\infty$.

Let $l_n=[\tau_n,\tau_{n+1}]$, and let $l$ be the union of $l_n$ for all $n$ such that $|\tau_n-\tau^{(0)}|<\delta$. Since the interval $[\tau^{(0)}-\delta,\tau^{(0)}+\delta]$ of integration in $I$ differs from $l$ by two small intervals near points $\tau^{(0)}\pm\delta$, and $S-\varphi>\gamma_2>0$ there, we have
\[
I=|x|^{(2-d)/2}\int_{l}h(\tau)e^{-|x|S(\tau,\theta,\lambda)}d\tau+O(e^{-|x|[\varphi(\theta,\lambda)+\gamma_2]}), \quad  |x|\to\infty.
\]
We represent the integral over $l$ as the sum of the integrals over $l_n$, and then write each of the latter integrals as the value of the integrand at an intermediate point $\tau_n+\sigma_n,~0\leq \sigma_n\leq 1/|x|,$ multiplied by $|l_n|=1/|x|$. This implies
\[
I=\sum_{n:~|\tau_n-\tau^{(0)}|<\delta}|x|^{-d/2}h(\tau_n+\sigma_n)e^{-|x|S(\tau_n+\sigma_n,\theta,\lambda)}+O(e^{-|x|[\varphi(\theta,\lambda)+\gamma_2]}),  \quad  |x|\to\infty.
\]

Since
$$
|S(\tau_n+\sigma_n,\theta,\lambda)-S(\tau_n,\theta,\lambda)|\leq \max_{|\tau-\tau^{(0)}|<\delta}|S_\tau|\sigma_n\leq \varepsilon/|x|,
$$
and $\varepsilon $ is small, we have
$$
e^{-|x|S(\tau_n+\sigma_n,\theta,\lambda)}=e^{-|x|S(\tau_n,\theta,\lambda)}(1+\eta_n), \quad |\eta_n|<2\varepsilon.
$$
The function $h$ is smooth and bounded from above and below by positive constants. Hence $h(\tau_n+\sigma_n)=h(\tau_n)(1+O(1/|x|))$. The last two relations allow us to rewrite $I$ in the form
\[
I=\sum_{n:~|\tau_n-\tau^{(0)}|<\delta}|x|^{-d/2}h(\tau_n)e^{-|x|S(\tau_n,\theta,\lambda)}(1+\xi_n)+O(e^{-|x|[\varphi(\theta,\lambda)+\gamma_2]}),  \quad  |\xi_n|<3\varepsilon,  \quad  |x|\to\infty.
\]
Hence
\begin{equation}\label{xxy}
\frac{\widehat{I}}{1+3\varepsilon}\leq \sum_{n:~|\tau_n-\tau^{(0)}|<\delta}|x|^{-d/2}h(\tau_n)e^{-|x|S(\tau_n,\theta,\lambda)}\leq \widehat{I}(1+3\varepsilon),   \quad  |x|\to\infty,
\end{equation}
where $\widehat{I}=I+O(e^{-|x|[\varphi(\theta,\lambda)+\gamma_2]})$. To complete the proof of the theorem, it remains only to replace $I$ by its asymptotics as $|x|\to\infty$ given by the Laplace method, and then combine (\ref{xxy}) (where $\varepsilon$ is arbitrary small) with (\ref{1222}).

\qed

\section{Propagation of the front}
Let $u=u(t,x)$ be the density of a population with initial density $u_0(x)$ being a non-negative continuous function with a compact support. It satisfies the equations
\begin{equation}\label{u0ca}
\frac{\partial u}{\partial t}=(\mathcal L+v(x))u, \quad \quad u(0,x)=u_0(x).
\end{equation}
We assume that the potential $v$ is non-negative, continuous and has a compact support. If the Hamiltonian $H=\mathcal L+v$ has a positive ground state energy, then the expected total population $\int_{R^d}u(t,x)dx$ grows with time, and an important question is to describe the location of the propagating front $F(t):=\{u(t,x):u=1\}$. This problem will be studied in this section when the transition density $a(y)$ has an ultra light tail, i.e., (\ref{ult}) holds. However, in the supplementary lemmas below we assume only that
the transition density $a(y)$ has a light tail, i.e., $|a(y)|<Ce^{-\delta|y|}$ with some $\delta>0$, and therefore $\widehat{a}(k)$ is analytic in $k$ when $|{\rm Im}k|<\delta$. We also will extend assumption (\ref{cond}) somewhat by requiring that
\begin{equation}\label{l11}
\|\widehat{a}(k+i\tau)\|_{L^1(R^d_k)}\leq C, \quad  |\tau|\leq\varepsilon_0,
\end {equation}
for some $\varepsilon_0< \delta.$

The following lemma can be proved by repeating the standard arguments used to show that Fourier transforms of $L_1$-functions decay at infinity.
\begin{lemma}\label{l61} Let $a(y)$ have a light tail. Then the
function $|\widehat{a}(z)|,~z=k+i\tau,$ decays uniformly in $\tau$ in each region $ |\tau|\leq \delta_1<\delta$ when $|k|\to\infty.$
\end{lemma}
The next lemma provides an estimate for the transition density $p$ given by (\ref{trd}).
\begin{lemma}\label{l62} Let $a(y)$ have a light tail and (\ref{l11}) hold. Then
there exist constants $\alpha<\infty$ and $C<\infty$ such that the following estimate holds
\[
p(t,x)\leq C te^{\alpha\varepsilon^2t-\varepsilon|x|}, \quad  x\neq 0,
\]
for small enough $\varepsilon\geq 0$.
\end{lemma}
{\bf Proof.} We write $p$ in the form (\ref{pmd}) and omit the first term on the right since we assume that $x\neq 0$. From Lemma \ref{l61} it follows that $R^d$ in the second term in the right-hand side of (\ref{pmd}) can be replaced by $R^d+i\tau, |\tau|<\delta.$ Thus
\begin{equation}\label{ppp}
p(t,x)=\frac{e^{-t}}{(2\pi)^d}\int_{R^d}w(t,z)e^{i(z,x)}dk, \quad z= k+i\tau, ~~|\tau|<\delta, \quad w:=e^{t\widehat{a}(z)}-1,~~x\neq 0.
\end{equation}

Note that
\[
|w(t,z)|\leq t|\widehat{a}(z)|, \quad z=k+i\tau,
\]
if $t|\widehat{a}(z)|\leq 1$ or $t|\widehat{a}(z)|> 1$ and Re$\widehat{a}(z)\leq 0$. If $t|\widehat{a}(z)|> 1$ and Re$\widehat{a}(z)>0$, then
\begin{equation}\label{lll}
|w(t,z)|\leq 2e^{t{\rm Re}\widehat{a}(z)}\leq 2t|\widehat{a}(z)|e^{t{\rm Re}\widehat{a}(z)}\leq Ct|\widehat{a}(z)|e^{t\sup{\rm Re}\widehat{a}(z)},  \quad z=k+i\tau,
\end{equation}
where the supremum is taken over $k\in R^d$. Since $\widehat{a}(0)=1$ (see (\ref{cond1})), we have that $\sup{\rm Re}\widehat{a}(z)>1/2$ when $|\tau|$ is small enough. Then (\ref{lll}) is valid in all the cases.
We will take $\tau=\varepsilon x/|x|,~ \varepsilon\leq \varepsilon_0$  with small enough $\varepsilon_0>0$. Then (\ref{ppp}), (\ref{lll}) and (\ref{l11}) imply that
\[
p\leq Cte^{t(\sup{\rm Re}\widehat{a}(z)-1)-\varepsilon|x|}, \quad z=k+i\tau,~~\varepsilon\leq \varepsilon_0.
\]

In order to complete the proof of Lemma \ref{l62}, it remains to show the existence of an $\alpha$ such that
\begin{equation}\label{ese}
{\rm Re}\widehat{a}(k+i\tau)-1\leq \alpha\varepsilon^2, \quad k\in R^d,
\end{equation}
if $\varepsilon\geq 0$ is small enough. Due to Lemma \ref{l61}, there exists $R<\infty$ such that (\ref{ese}) holds with $\alpha=0$ when $|k|\geq R$. Estimate (\ref{ese}) with $ |k|\leq R$ and $\varepsilon\geq 0$ is small enough follows from the Taylor expansion since $\widehat{a}(z)$ is smooth (analytic), $\widehat{a}(k)$ is real, and Re$\widehat{a}(k)\leq 1$.

\qed

The asymptotics of the Green function $G_\lambda(x),\lambda>0,$ as $ |x|\to\infty$, will be used below (see Theorem \ref{tgr}) and therefore, it will be assumed that $a(y)$ has an ultra light tail, i.e., (\ref{ult}) holds.

\begin{theorem}\label{tla}
Let $u$ be the solution of the problem (\ref{u0ca}),
where the transition density $a(y)$ has an ultra light tail, (\ref{l11}) holds, and operator $H=\mathcal L+v(x)$ has a unique positive eigenvalue $\lambda=\lambda_0$.

Then the front $F(t)$ has the form $|x|=\frac{\lambda_0t+\frac{1-d}{2}\ln t}{\varphi(\theta,\lambda_0)}+O(1), ~t\to\infty$, where $\varphi$ is defined in (\ref{fi}). The density $u$ grows (decays) exponentially in time uniformly in $x$ in any region inside  (outside, respectively) the front whose distance from the front exceeds $\gamma t$ with some $\gamma>0$.
\end{theorem}
{\bf Remark.} In the case of light tails, the total number of positive eigenvalues of operator $H$ is at most finite. This will be proved elsewhere. If there are several positive eigenvalues, then each of them generates its own front. Since it is not clear which of them propagates faster, the total front of the population can be obtained as the maximum of the fronts generated by individual eigenvalues.

{\bf Proof of Theorem \ref{tla}.} The spectrum of the operator $H=\mathcal L+v$ consists of the interval $[-a,0]$ and a point $\lambda_0>0$. Using the spectral decomposition theorem, we obtain that
\begin{equation}\label{ab}
u(t,x)=\int_{-a}^0e^{\lambda t}dE_\lambda u_0 d\lambda+c_0\phi_0(x)e^{\lambda_0t},
\end{equation}
where $E_\lambda$ is the spectral projection operator, $\phi_0$ is a normalized eigenfunction with the eigenvalue $\lambda_0$, and $c_0=(u_0,\phi_0)$. Denote by $w$ the first term in the right-hand side above. Obviously, $\|w\|_{L^2}\leq \|u_0\|_{L^2}$ for all $t\geq 0$. We will show that the front is defined mostly by the second term in (\ref{ab}). An estimate on $w$ is needed to justify this fact.

Function $w$ satisfies the relations
\[
\frac{\partial w}{\partial t}-\mathcal Lw=-v(x)w, \quad \quad u(0,x)=u_0(x)-c_0\phi_0(x),
\]
and therefore $w=u_1+u_2$, where
\[
u_1=\int_0^t\int_{R^d}p(t-\tau, x-y)v(y)w(\tau,y)dyd\tau, \quad u_2=\int_{R^d}p(t, x-y)[u_0(y)-c_0\phi_0(y)]dy.
\]

Since $\|vw\|_{L^1}\leq C\|w\|_{L^2}\leq C$, from Lemma \ref{l62} it follows that
\begin{equation}\label{u1}
|u_1|\leq C\int_0^t te^{\alpha\varepsilon^2t-\varepsilon|x|}d\tau=Ct^2e^{\alpha\varepsilon^2t-\varepsilon|x|}
\end{equation}
when $x$ does not belong to the support of $v$ and $\varepsilon\geq 0$ is small enough. We write $u_2$ as the sum $u_{21}+u_{22}$ by separating the terms in the square brackets in the expression for $u_2$. Obviously, an estimate similar to (\ref{u1}) (with the factor $t^2$ replaced by $t$) is valid for $u_{21}$ when $x $  does not belong to the support of $u_0$ and $\varepsilon\geq 0$ is small enough. In order to estimate $u_{22}$, we note that equation $(\mathcal L+v)\phi_0=-\lambda_0 \phi_0$ implies that $\phi_0=-G_{\lambda_0}(v\phi_0)$ has the same asymptotics at infinity as $G_{\lambda_0}$, which is given by Theorem \ref{tgr}, i.e.,
\begin{equation}\label{phi0}
\phi_0= g(\theta, \lambda_0)|x|^{(1-d)/2}e^{-|x|\varphi(\theta,\lambda_0)}(1+o(1)),  \quad |x|\to\infty.
\end{equation}
From (\ref{pmd}), Lemma \ref{l62}, and (\ref{phi0}), it follows that
\[
|u_{22}+c_0e^{-t}\phi_0(x)|\leq  C\int_{R^d} te^{\alpha\varepsilon^2t-\varepsilon|x-y|}e^{-\gamma |y|}dy
\]
with some $\gamma>0$. Thus,
\[
|u_{22}|\leq  C(1+t)e^{\alpha\varepsilon^2t-\varepsilon|x|}.
\]
By combining estimates for $u_1,u_{21},u_{22}$, we obtain that
\begin{equation}\label{u-}
|w|\leq  C(1+t^2)e^{\alpha\varepsilon^2t-\varepsilon|x|}
\end{equation}
when $x$ does not belong to the supports of $v$ and $u_0$, and $\varepsilon\geq 0$ is small enough.

In order to complete the proof of the theorem, we split $R^d$ in the following regions that depend on time. The central region $D_0\subset R^d$ is defined by  inequalities
\[
\frac{1}{2}<c_0\phi_0(x)e^{\lambda_0t}<2.
\]
Asymptotics (\ref{phi0}) implies that the solution of the equation $c_0\phi_0(x)e^{\lambda_0t}=c$ has the form $|x|=\frac{\lambda_0t+\frac{1-d}{2}\ln t}{\varphi(\theta,\lambda_0)}+O(1), ~t\to\infty$. We denote the remainder term by $O_1(1), O_2(1)$ when $c=1/2,2$, respectively. Here $O_1(1)> O_2(1)$, and $D_0$ can be rewritten, for large $t$, as follows
\[
\frac{\lambda_0t+\frac{1-d}{2}\ln t}{\varphi(\theta,\lambda_0)}+O_2(1)<|x|<\frac{\lambda_0t+\frac{1-d}{2}\ln t}{\varphi(\theta,\lambda_0)}+O_1(1), \quad t\to\infty.
\]
By $D_1\subset R^d$ we denote the region where one of the following inequalities hold as $t\to\infty$:
\[
\frac{(\lambda_0-\gamma)t}{\varphi(\theta,\lambda_0)}\leq|x|\leq \frac{\lambda_0t+\frac{1-d}{2}\ln t}{\varphi(\theta,\lambda_0)}+O_2(1) \quad {\rm or} \quad \frac{\lambda_0t+\frac{1-d}{2}\ln t}{\varphi(\theta,\lambda_0)}+O_1(1)\leq|x|\leq\frac{(\lambda_0+\gamma)t}{\varphi(\theta,\lambda_0)},
\]
where $\gamma$ is arbitrary small. Finally, let
\[
D_2=\{x:~ |x|<\frac{(\lambda_0-\gamma)t}{\varphi(\theta,\lambda_0)}\}, \quad D_3=\{x:~ |x|>\frac{(\lambda_0+\gamma)t}{\varphi(\theta,\lambda_0)}\}, \quad t\to\infty.
\]

Since $c_0\phi_0(x)e^{\lambda_0t}$ grows exponentially in $t$ when $x\in D_2$, estimate (\ref{u-}) with $\varepsilon=0$ implies that $u,~ x\in D_2,$ grows exponentially as $t\to\infty$. Since $c_0\phi_0(x)e^{\lambda_0t}$ decays exponentially in $t$ when $x\in D_3$, estimate (\ref{u-}) with small enough $\varepsilon>0$ implies that $u,~ x\in D_3,$ decays exponentially as $t\to\infty$. Thus front $F(t)$ belongs to $D_1\bigcup D_0$. From (\ref{u-}) with small enough $\varepsilon=0$ it follows that $w$ decays exponentially when $x\in D_1,~t\to \infty$. Hence $u>\frac{3}{2}$ or $u<\frac{3}{4}$ when $x\in D_1$ and $t$ is large enough, i.e., $F(t)\subset D_0$ as $t\to \infty$.

\qed
\section{Limit of the population density as $t\to \infty$.}
Assume now that the Hamiltonian $H=\mathcal L+w$ does not have positive spectrum. Then one can expect that the total population does not grow with time and the population density has a limit as $t\to\infty$. We will prove the latter facts under a somewhat stronger assumption than the absence of the positive spectrum. Recall (see Theorem \ref{sss}) that the absence of the positive spectrum requires the underlying process with the generator $\mathcal L$ to be transient. Moreover, the potential has to be small enough. Additionally, we will assume that the norm of the operator $G_0v$ in the space $C=C_b(R^d)$ of continuous bounded functions ($\|u\|_C={\rm sup}_x|u|$) is less than one. Here $G_0$ is the operator defined in (\ref{gg}) and $G_0v\psi:=G_0(v\psi)$. At the end of the section, we will provide an estimate on $v$ that implies  $\|G_0v\|<1$.
We do not impose any restrictions on the tail of the transition density $a(y)$.

We will consider the population with a constant initial distribution, i.e., the density $u=u(t,x)\in C'=\bigcap_{T>0}C^T$ is the solution of the problem
\begin{equation}\label{u0c1}
\frac{\partial u}{\partial t}=(\mathcal L+v(x))u, \quad \quad u(0,x)=1.
\end{equation}
Here $T$ is arbitrary and $C^T$ is the space of functions that are continuous in $(t,x)$ and bounded when $x\in R^d,0\leq t\leq T$, with $\|u\|_{C^T}={\rm sup}_{(t,x)} |u|$. We are not going to discuss the uniqueness in the space $C'$. Instead we study a specific solution in $ C'$ that can be obtained from the corresponding integral equation (see more details in the proof of the theorem below).
\begin{theorem}\label{tlast}
Let $H=\mathcal L+v$ not have positive spectrum and $\|G_0v\|_C<1$. Then the solution $u=u(t,x)\in C'$ of problem (\ref{u0c1}) has a limit in the space $C'$ as $t\to\infty$ and $\lim_{t\to\infty}u=(1-G_0v)^{-1}1$.
\end{theorem}
We will need the following analogue of Lemma \ref{l62}. Denote by $P_tf$ solution $u=u(t,x)$ of the unperturbed problem
\[
u_t=\mathcal Lu, \quad t>0, \quad u(0,x)=f(x)\in C_{\rm com},
\]
that is obtained by convolution with the fundamental solution (\ref{trd}):
\[
u=P_tf=\int_{R^d} p(t,x-y)f(y)dy.
\]
\begin{lemma}\label{llast}
Let $p(t,x)$ be the transition density of a transient random walk with the generator $\mathcal L$. Then for each $f\in C $ with the support in the ball $|x|<R$ and each $T<\infty$, the function $P_tf$ belongs to $C^T$, and
\[
\|P_tf\|_{C^T}\leq C(R,T)\|f\|_C.
\]
\end{lemma}
{\bf Proof.} From (\ref{pmd}) it follows that $P_tf=e^{-t}f+e^{-t}F^{-1}w$, where the second term is the inverse Fourier transform of $w=[e^{t\widehat{a}(k)}-1]\widehat{f}(k)$. It is enough to prove the statement of the lemma for the term $F^{-1}w$. Since $|\widehat{f}(k)|\leq C(R)\|f\|_C$ and $\|e^{t\widehat{a}(k)}-1\|_{L^1}\leq C\|T\widehat{a}(k)\|_{L^1}=CT$ (due to (\ref{cond})), it follows that $F^{-1}w$ is continuous in $x$ and $|F^{-1}w|\leq C(R,T),~t\leq T.$ Obviously, $w$ as element of $L^1$, depends smoothly on $t$. This implies that $F^{-1}w$ is continuous in $(t,x)$.

\qed

{\bf Proof of Theorem \ref{tlast}.} From the Duhamel's principle, it follows that every solution $u\in C'$ of (\ref{u0c1}) satisfies the equation
\begin{equation}\label{vv}
u=1+\int_0^tP_{t-s}(vu(s,\cdot))ds.
\end{equation}
Equation (\ref{vv}) is uniquely solvable in $C^T$, and its solution satisfies (\ref{u0c1}). Theorem \ref{tlast} concerns this solution $u$.
Indeed, the unique solvability of (\ref{vv}) follows immediately from the facts that equation (\ref{vv}) is of the Volterra type and the operator $P_tv:C\to C$ is bounded uniformly in $t\in[0,T]$. Moreover, solution  $u\in C^T$ can be obtained by iterations:
\begin{equation}\label{uq}
u=\sum _0^\infty (Q_t)^n1, \quad Q_tf=\int_0^tP_{t-s}(vf(s,\cdot), \quad t\leq T.
\end{equation}
It remains to find the limit of $u$ as $t\to\infty$.

Assume that $f(s,x)$ has the following properties: $f\geq 0,~f(s_1,x)\geq f(s_2,x)$ when $s_1>s_2$, and there is a limit $f_0(x):=\lim_{s\to\infty}f(s,x)\in C$. Then $Q_tf$ has the same properties as $f$ and
\begin{equation}\label{up}
Q_tf\uparrow G_0(vf_0) \quad {\rm as} \quad t\to \infty.
\end{equation}
Indeed,
\begin{equation}\label{mon}
Q_tf=\int_0^t p(s,x-y)v(y)f(t-s,y)dyds.
\end{equation}
The latter relation and Lemma \ref{l1} imply that $Q_tf\geq 0$ and $Q_tf$ is monotone in $t$. An upper bound for $Q_tf$ can be obtained by putting $t=\infty$ in the right-hand side of (\ref{mon}). Since the process with the transition density $p(t,x)$ is transient, it follows that
$$
\int_0^\infty p(s,x-y)ds=G_0(x-y),
$$
and therefore $Q_tf\leq G_0(vf_0).$  Thus $Q_tf$ has a limit as $t\to\infty$ and $\lim_{t\to\infty}Q_tf\leq G_0(vf_0).$ The lower bound for $\lim_{t\to\infty}Q_tf$ is the same, as also follows from  (\ref{mon}):
\[
Q_tf\geq \int_0^{t/2} p(s,x-y)v(y)f(t-s,y)dyds\geq \int_0^{t/2} p(s,x-y)v(y)f(t/2,y)dyds.
\]
By passing to the limit in the last inequality as $t\to\infty$, we obtain that $\lim_{t\to\infty}Q_tf\geq G_0(vf_0).$ Hence, (\ref{up}) is proved.

Using (\ref{up}), one can easily show by induction that $(Q_t)^n1\uparrow (G_0v)^n1$. This and (\ref{uq}) complete the proof of the theorem.

\qed

We will conclude this section by providing a sufficient condition for the estimate $\|G_0v\|<1.$
\begin{lemma}
The following estimate holds
\[
\|G_0v\|_C\leq \|v\|_C+\frac{1}{(2\pi)^d}\int_{R^d}\frac{|\widetilde{a}(k)\widetilde{v}(k)|}{1-a(k)}dk.
\]
\end{lemma}
The statement follows immediately from the relation $G_0=I+T_0$, where $T_0$ is given by (\ref{tl}).

\end{document}